\documentstyle[auxdefs,diagram,12pt]{article}

\begin{document}

\input amssym.def

\input amssym

\newcommand{\ba}{\mbox{$\Bbb A$}}
\newcommand{\bz}{\mbox{$\Bbb Z$}}
\newcommand{\bq}{\mbox{$\Bbb Q$}}
\newcommand{\bc}{\mbox{$\Bbb C$}}
\newcommand{\bh}{\mbox{$\Bbb H$}}
\newcommand{\br}{\mbox{$\Bbb R$}}
\newcommand{\bp}{\mbox{$\Bbb P$}}
\newcommand{\bres}{\mbox{$\Bbb F$}}

\newcommand{\ra}{\mbox{$\rightarrow$}}
\newcommand{\lra}{\mbox{$\longrightarrow$}}
\newcommand{\proof}{\noindent {\it Proof: }}
\newcommand{\remark}{\noindent {\bf Remark: }}

\newcommand{\cA}{\mbox{${\cal A}$}}
\newcommand{\cB}{\mbox{${\cal B}$}}
\newcommand{\cC}{\mbox{${\cal C}$}}
\newcommand{\cD}{\mbox{${\cal D}$}}
\newcommand{\cE}{\mbox{${\cal E}$}}
\newcommand{\cF}{\mbox{${\cal F}$}}
\newcommand{\cG}{\mbox{${\cal G}$}}
\newcommand{\cH}{\mbox{${\cal H}$}}
\newcommand{\cI}{\mbox{${\cal I}$}}
\newcommand{\cJ}{\mbox{${\cal J}$}}
\newcommand{\cK}{\mbox{${\cal K}$}}
\newcommand{\cL}{\mbox{${\cal L}$}}
\newcommand{\cM}{\mbox{${\cal M}$}}
\newcommand{\cN}{\mbox{${\cal N}$}}
\newcommand{\cO}{\mbox{${\cal O}$}}
\newcommand{\cP}{\mbox{${\cal P}$}}
\newcommand{\cQ}{\mbox{${\cal Q}$}}
\newcommand{\cR}{\mbox{${\cal R}$}}
\newcommand{\cS}{\mbox{${\cal S}$}}
\newcommand{\cT}{\mbox{${\cal T}$}}
\newcommand{\cU}{\mbox{${\cal U}$}}
\newcommand{\cV}{\mbox{${\cal V}$}}
\newcommand{\cW}{\mbox{${\cal W}$}}
\newcommand{\cX}{\mbox{${\cal X}$}}
\newcommand{\cY}{\mbox{${\cal Y}$}}
\newcommand{\cZ}{\mbox{${\cal Z}$}}

\setlength{\unitlength}{1cm} \thicklines \newcommand{\ci}{\circle*{.13}}
\renewcommand{\mp}{\multiput}
\newcommand{\num}[1]{{\raisebox{-.5\unitlength}{\makebox(0,0)[b]{${#1}$}}}}
\newcommand{\mynum}[1]{{\raisebox{-.1\unitlength}{\makebox(-.4,0)[r]{$Y_{#1}$}}}}

\newcommand{\da}{\downarrow}
\newcommand{\el}{\ell}

\newcommand{\aut}{automorphism}
\newcommand{\clg}{compact Lie group}
\newcommand{\defr}{deformation retract}
\newcommand{\diffeo}{diffeomorphism}
\newcommand{\ems}{Eilenberg-Maclane}
\newcommand{\fg}{finitely-generated} 
\newcommand{\fd}{finite dimensional}

\newcommand{\homeo}{homeomorphism}
\newcommand{\ho}{homomorphism}
\newcommand{\homeq}{homotopy equivalence}
\newcommand{\hty}{homotopy type}
\newcommand{\iso}{isomorphism}
\newcommand{\les}{long exact sequence}
\newcommand{\lfp}{Lefschetz fixed-point theorem}
\newcommand{\mvs}{Mayer-Vietoris sequence}
\newcommand{\ses}{short exact sequence}
\newcommand{\nbhd}{neighborhood}
\newcommand{\pid}{principal ideal domain}
\newcommand{\sss}{spectral sequence}
\newcommand{\we}{weak equivalence}

\newcommand{\wrt}{with respect to}

\newcommand{\rn}{\mbox{$\br ^n$}}
\newcommand{\cn}{\mbox{$\bc ^n$}}

\newcommand{\gc}{\mbox{$G_{\Bbb C}$}}

\newcommand{\phiplus}{\mbox{$\Phi ^+$}}
\newcommand{\phiminus}{\mbox{$\Phi ^-$}}
\newcommand{\phii}{\mbox{$\Phi _I$}}
\newcommand{\phij}{\mbox{$\Phi _J$}}
\newcommand{\phiiplus}{\mbox{$\Phi _I ^+$}}
\newcommand{\phiiminus}{\mbox{$\Phi _I ^-$}}
\newcommand{\phijplus}{\mbox{$\Phi _J ^+$}}
\newcommand{\phijminus}{\mbox{$\Phi _J ^-$}}
\newcommand{\phiirad}{\mbox{$\Phi _I ^{rad}$}}
\newcommand{\phijrad}{\mbox{$\Phi _J ^{rad}$}}
\newcommand{\phiipar}{\mbox{$\Phi _I ^{par}$}}
\newcommand{\phijpar}{\mbox{$\Phi _J ^{par}$}}
\newcommand{\ujrad}{\mbox{$U _J ^{rad}$}}
\newcommand{\uirad}{\mbox{$U _I ^{rad}$}}

\newcommand{\lieuirad}{u_I ^{rad}}
\newcommand{\lieujrad}{u_J ^{rad}}
\newcommand{\lieualpha}{u_{\alpha}}
\newcommand{\lieqw}{q_w}
\newcommand{\liepi}{p_I}
\newcommand{\lieuw}{u_w}
\newcommand{\lieuwp}{u_w ^\prime}
\newcommand{\liegc}{\frak{g}_{\Bbb C}}
\newcommand{\liehc}{\frak{h}_{\Bbb C}} 
\newcommand{\gctil}{\mbox{$\tilde{G} _{\Bbb C}$}}
\newcommand{\gctilp}{\mbox{$\tilde{G} _{\Bbb C}/P$}}
\newcommand{\gctilad}{\mbox{$\tilde{G} _{\Bbb C}^{ad}$}}
\newcommand{\gctiladp}{\mbox{$\tilde{G} _{\Bbb C} ^{ad}/P^{ad}$}}
\newcommand{\btil}{\mbox{$\tilde{B}$}}
\newcommand{\btilminus}{\mbox{$\tilde{B}^-$}}  
\newcommand{\dtil}{\mbox{$\tilde{D}$}}
\newcommand{\util}{\mbox{$\tilde{U}$}}
\newcommand{\utilminus}{\mbox{$\tilde{U} ^-$}}
\newcommand{\utilminussig}{\mbox{$\tilde{U} ^- _\sigma$}}
\newcommand{\wtil}{\mbox{$\tilde{W}$}}
\newcommand{\wtils}{\mbox{$\tilde{W}^S$}}
\newcommand{\ntil}{\mbox{$\tilde{N}$}}
\newcommand{\stil}{\mbox{$\tilde{S}$}}
\newcommand{\liegctil}{\mbox{$\tilde{\frak{g}} _{\bc}$}}
\newcommand{\phitil}{\mbox{$\tilde{\Phi}$}}
\newcommand{\phitilplus}{\mbox{$\tilde{\Phi} ^+$}}
\newcommand{\phitilminus}{\mbox{$\tilde{\Phi} ^-$}}
\newcommand{\coroot}{\mbox{$\cQ ^\vee$}}
\newcommand{\antidom}{\mbox{$Q^\vee _-$}}
\newcommand{\chamber}{\overline{\cC}}
\newcommand{\affg}{\mbox{$\cL _G$}}   
\newcommand{\tchat}{\mbox{$\hat{T} _{\bc}$}}
\newcommand{\that}{\mbox{$\hat{T}$}}  

\newcommand{\rvar}{Richardson variety}
\newcommand{\rvars}{Richardson varieties}  
\newcommand{\rfilt}{Richardson filtration} 
\newcommand{\bvar}{Birkhoff variety}
\newcommand{\bvars}{Birkhoff varieties}
\newcommand{\cpo}{closed parabolic orbit}
\newcommand{\clr}{coroot lattice representative}  
\newcommand{\ef}{equivariant formality}
\newcommand{\svar}{Schubert variety}
\newcommand{\svars}{Schubert varieties}
\newcommand{\pd}{Poincar\' e duality}
\newcommand{\pds}{Poincar\' e duality space}
\newcommand{\ppoly}{Poincar\' e polynomial}
\newcommand{\pser}{Poincar\' e series}
\newcommand{\psg}{one-parameter subgroup}
\newcommand{\fnk}{\mbox{$X_{n,k}$}}
\newcommand{\fnkp}{\mbox{$X_{n,k} ^\flat$}}
\newcommand{\ppc}{positive pair condition}
\newcommand{\npc}{negative pair condition}
\newcommand{\osc}{opposite sign condition}
\newcommand{\und}{\underline{0}}
\newcommand{\nope}{$\lambda \notin \cQ ^\vee$}
\newcommand{\onemin}{1- \alpha _0 (\lambda)} 
\newcommand{\owt}{(otherwise $\lambda$ is overweight)} 
\newcommand{\adp}{admissible palindromic}
\newcommand{\fork}{(otherwise $\lambda$ covers a fork)}

\newcommand{\gsm}{graph-splitting move} 
\newcommand{\aslam}{\alpha _s
(\lambda)} 
\newcommand{\azlam}{\alpha _0 (\lambda)}
\newcommand{\atlam}{\alpha _t (\lambda)} 
\newcommand{\aulam}{\alpha _u
(\lambda)} 
\newcommand{\astlam}{\alpha _s (\lambda) + \alpha _t
(\lambda)} 
\newcommand{\alam}{\alpha (\lambda)} 
\newcommand{\blam}{\beta
(\lambda)} 
\newcommand{\aplam}{\alpha ^\prime (\lambda)}
\newcommand{\bpp}{$\beta$-positive pair}
\newcommand{\bnp}{$\beta$-negative pair}
\newcommand{\lam}{\mbox{$\lambda$}}
\newcommand{\circular}{spiral}

\newcommand{\xlam}{\mbox{$X_\lambda$}}
\newcommand{\zlam}{\mbox{$Z_\lambda$}}
\newcommand{\elam}{\mbox{$e_\lambda$}}
\newcommand{\slam}{\mbox{$S_\lambda$}}
\newcommand{\openlam}{\mbox{$\cU _\lambda$}}
\newcommand{\openslam}{\mbox{$\cS _\lambda$}}
\newcommand{\nlam}{\mbox{$N_\lambda$}}
\newcommand{\xmu}{\mbox{$X_\mu$}}
\newcommand{\zmu}{\mbox{$Z_\mu$}}
\newcommand{\emu}{\mbox{$e_\mu$}}
\newcommand{\smu}{\mbox{$S_\mu$}}
\newcommand{\openmu}{\mbox{$\cU _\mu$}}

\newcommand{\vlamu}{\mbox{$V_{\lambda ,\mu}$}}
\newcommand{\klamu}{\mbox{$K_{\lambda ,\mu}$}}
\newcommand{\vlamuplus}{\mbox{$V_{\lambda ,\mu} ^+$}}
\newcommand{\vlam}{\mbox{$\cE _\lambda$}}
\newcommand{\vi}{\mbox{$\cE _{\cal I}$}}
\newcommand{\zi}{\mbox{$Z _{\cal I}$}}
\newcommand{\zip}{\mbox{$Z _{\cal I} ^\prime$}}
\newcommand{\xj}{\mbox{$X _{\cal J}$}}
\newcommand{\sj}{\mbox{${\cal S} _{\cal J}$}}   

\newcommand{\utilam}{\mbox{$\tilde{U} _\lambda$}}
\newcommand{\utilamminus}{\mbox{$\tilde{U} _\lambda ^-$}}
\newcommand{\utilamp}{\mbox{$\tilde{U} _\lambda ^\prime$}} 
 \newcommand{\utilampminus}{\mbox{$\tilde{U} _\lambda ^{\prime -}$}} 
\newcommand{\lieutilam}{\mbox{$\frak{u} _\lambda ^\prime$}}  
\newcommand{\lieutilamp}{\mbox{$\tilde{U} _\lambda$}}  
\newcommand{\ulamk}{\mbox{$\tilde{U} _\lambda$}}  
\newcommand{\ulampk}{\mbox{$\tilde{U} _\lambda$}}  

\newcommand{\balp}{\mbox{$\gctil \times _P \affg$}} 
\newcommand{\esiglam}{\mbox{$E_{\sigma, \lambda}$}}  
\newcommand{\genc}{generating complex}
\newcommand{\ggenc}{geometric generating complex} 
\newcommand{\genv}{generating variety}  
\newcommand{\ggenv}{geometric generating variety}  
\newcommand{\genvs}{generating varieties}  
\newcommand{\xlamz}{\mbox{$X_{\lambda _0}$}}  

\newcommand{\LG}{\cL _G}
\newcommand{\C}{\Bbb C}   
\newcommand{\Z}{\Bbb Z}   
\newcommand{\GC}{G_{\Bbb C}}   
\newcommand{\uh}{\mbox{$\frak{u}_{\alpha _0}$}}
\newcommand{\G}{\tilde{G} _{\Bbb C}}

\newtheorem{theorem}{Theorem}[section]

\newtheorem{proposition}[theorem]{Proposition}

\newtheorem{lemma}[theorem]{Lemma}

\newtheorem{conjecture}[theorem]{Conjecture}

\newtheorem{example}[theorem]{Example}

\newtheorem{corollary}[theorem]{Corollary}

\newtheorem{definition}[theorem]{Definition}

\title{Generating varieties for affine Grassmannians}

\author{Peter J. Littig and Stephen A. Mitchell}

\maketitle

\section{Introduction} 

Let $G$ be a simply-connected compact Lie group with simple Lie algebra,
and let $G_{\C}$ be its complexification.  The affine Grassmannian
associated to $G$ is the homogeneous space $\LG = \G/P$, where $\G$ is
the affine Kac-Moody group of regular maps $\C^{\times} \lra G_{\C}$ and
$P$ is the subgroup of maps $\C \lra G_{\C}$.  The space $\LG$ is a
projective ind-variety and has the homotopy type of $\Omega G$, the
space of based loops on $G$. Sitting inside $\LG$ is a family of
finite-dimensional, complex projective varieties, called {\it affine
Schubert varieties}.  These varieties are indexed by the elements of
\coroot , the coroot lattice associated to $G$.  We write $X_{\lambda}$
for the variety corresponding to $\lambda \in \coroot$.

Its infinite-dimensionality notwithstanding, \affg\ is in many ways
analogous to an ordinary Grassmannian. There are, however, two striking
differences: (1) \affg\ is a topological group; and (2) the principal
bundle defining \affg\ as a homogeneous space is topologically
trivial. No finite-dimensional flag variety $\gc /Q$ has these
properties. These two properties of $\cL _G$ give rise to interesting
new phenomena that have no analogue in the classical setting. The most
obvious is that there is a Schubert calculus not only for the cup
product in cohomology, but also for the Pontrjagin product in
homology. The problem is to compute the coefficients $a_{\sigma, \eta}
^\mu$ in the expansion

$$[X_\sigma] \cdot [X_\eta] =\sum a_{\sigma, \eta }^ \mu [X_\mu].$$

As far as we know, this problem was first considered in the second
author's 1987 paper \cite{mbottfilt} on minuscule generating complexes,
albeit in a very special case with specific topological applications in
mind. Since then it has been studied by several authors, notably in
upublished work of Dale Peterson and in \cite{llms}, \cite{lss},
\cite{magyar}.

In the present paper we consider the Pontrjagin product on the
point-set, geometric level. A finite-dimensional projective subvariety
$X \subset \affg$ is a {\it generating variety} if the image of $H_* X
\lra H_* \affg$ generates $H_* \affg$ as a ring. A \svar\ \xlam\ is a {\it
geometric generating variety} if (1) for all $n\geq 1$ the image of the
$n$-fold multiplication $X_\lambda ^n \lra \affg$ is a
\svar\ $X_{\lambda _n}$, with $dim \, X_{\lambda _n} =n \cdot dim \,
\xlam$; and (2) $H_* X_\lambda ^n \lra H_* X_{\lambda _n}$ is
surjective. Let $\alpha _0$ denote the highest root of $G$. 

\begin{theorem} \label{intro} \marginpar{}

$X_{-\alpha _0 ^\vee}$ is a \ggenv , with $\lambda _n =-n \alpha _0
  ^\vee$. 

\end{theorem}

A classical theorem of Bott \cite{bott} can be reinterpreted to say that
every \affg\ admits a smooth generating variety. In contrast, the
\genvs\ of the theorem above are all singular. Their advantage is that
they are \svars , and give one canonical \genv\ in each Lie type. It
turns out that \affg\ admits a smooth Schubert \genv\ if and only if $G$
is not of type $E_8$, $F_4$, or $G_2$ (Proposition~\ref{notsmooth}).
In type $C_n$, it turns out that $X_{-\alpha _0 ^\vee}$ is essentially
the generating complex for $\Omega Sp(n)$ considered by Hopkins
\cite{hopkins}. 

The proof of Theorem~\ref{intro} makes use of the balanced product
$\gctil \times _P \affg$. We construct an ind-variety structure on the
balanced product, whose filtrations are topologically trivial algebraic
fiber bundles $E_{\sigma ,\lambda}$ with base $X_\sigma$ and fiber
$X_\lambda$ (\lam\ an anti-dominant element of the coroot lattice). Then
we relate the group multiplication $\affg \times \affg \lra \affg$ to
the action map $\gctil \times _P \affg \lra \affg$.

Along the way we consider the surprising topological properties (1) and
(2) and how they fail in the category of ind-varieties. In particular,
we show that there is no topological group structure on \affg\ whose
multiplication map is a morphism of ind-varieties
(Corollary~\ref{nonexist}). Regarding property (2), although it is clear
that the principal $P$-bundle $\gctil \lra \gctil /P$ is not trivial as
ind-varieties, it is also easy to show that the balanced product bundle
is isomorphic to $\affg \times \affg$ as ind-varieties. It does not
follow, however, that the bundles $E_{\sigma ,\lambda}$ are products. 

\section{Notation and Conventions}

Throughout this paper, $G$ is a simple, simply-connected, compact Lie
group with maximal torus $T$, Weyl group $W$, and complexification
$G_{\C}$.  We fix a Borel group $B$ containing the Cartan subgroup $H =
T_{\C}$, and let $B^-$ denote the opposite Borel group.  Let $U$ be the
unipotent radical of $B$, so that $B = HU$. Let $S \subset W$ denote the
simple reflections, root system $\Phi$, and simple roots $\alpha _s,
s\in S$. Let \coroot\ denote the coroot lattice.  We denote by
$\mathcal{C}$ the dominant Weyl chamber determined by $B$.

\bigskip
\noindent
{\it Loop Groups.} The primary references for this section
are \cite{ps} and \cite{kumar}; see also \cite{mbuilding}. 

An {\it algebraic loop in} $G_{\C}$ is a regular map $\C^{\times} \lra
G_{\C}$.  The set of all such loops, which we denote by $\G$, is a
topological group with respect to pointwise multiplication.  In type
$A_{n-1}$, for example, $G_{\C} = SL_n \C$ and $\G =
SL_n(\C[z,z^{-1}])$.  In general Lie type, $\G$ can be shown to be a
complex affine ind-group.

Let $\Omega_{alg} G$ be the subgroup
$$
\Omega_{alg} G = \{ f \in \widetilde{G}_{\C} : f|_{S^1} \in \Omega G
\}.
$$

\noindent
We call $\Omega_{alg} G$ the group of based algebraic loops on $G$.
Now, let $P$ be the subgroup of $\widetilde{G}_{\C}$ of algebraic loops $\C
\lra G_{\C}$ (i.e., algebraic loops that may be extended over the
origin.)  The following decomposition theorem is a fundamental part of
the structure theory of loop groups.

\begin{theorem} \label{Iwasawa}
[Iwasawa Decomposition] The multiplication map $\Omega_{alg} G \times
  P \lra \G$ is a homeomorphism.
\end{theorem}

We call this the Iwasawa decomposition because it is a modified analogue
of Iwasawa's KAN decomposition for real Lie groups. The unmodified
analogue takes the form $L_{alg} G \times A \times \util
\stackrel{\cong}{\lra} \gctil$, where $A$ comes from the $KAN$
decomposition of the constant loops \gctil , and $\util =\{ f\in \btil :
f(0) \in U^-\}$. 

Two corollaries of Theorem \ref{Iwasawa} are particularly important
for our purposes.  First, the homogeneous space $\affg$ is homeomorphic
to $\Omega_{alg} G$ and hence carries the structure of a topological
group.  As we shall observe later, such a group structure is wholly
absent in the classical setting of homogeneous spaces of complex
algebraic groups.  Second, the principal bundle $P \lra \gctil \lra \affg$ is
trivial.

The following theorem, due, independently, to Quillen and
Garland-Raghunathan, establishes the connection between $\Omega_{alg}
G$ and the ordinary loop group $\Omega G$.

\begin{theorem} \label{Quillen}
[Quillen, Garland-Raghunathan] The inclusion $\Omega_{alg} G \lra
\Omega G$ is a homotopy equivalence.
\end{theorem}

Next, we describe the Tits system associated to $\G$.  Let
$\widetilde{B} = \{ f \in P : f(0) \in B^- \}$.  The group
$\widetilde{B}$ is called the {\it Iwahori subgroup} of
$\widetilde{G}_{\C}$ and may be thought of as an affine analog of a
Borel subgroup of $G_{\C}$.  Let $\widetilde{N}_{\C}$ denote the group
of algebraic loops in $G_{\C}$ that take their image in $N_{\C}$.  The
intersection $\widetilde{B} \cap \widetilde{N}_{\C}$ is normal in
$\widetilde{N}_{\C}$ and the quotient group $\widetilde{W} =
\widetilde{N}_{\C} / (\widetilde{B} \cap \widetilde{N}_{\C})$ is
canonically identified with the affine Weyl group associated to
$G$. Thus \wtil\ has a semi-direct product decomposition $\wtil =
\coroot \rtimes W$. Here we note that since $G$ is simply-connected,
the lattice $\cQ ^\vee \subset \frak{t} _{\Bbb C}$ can be identified
with the subgroup $Hom \, (S^1, T) \subset \Omega _{alg}
G$. Consequently, we alternate frequently between additive notation and
multiplicative notation; the meaning will be clear from the context.  It
is shown in \cite{mbuilding} that $(\widetilde{G}_{\C}, \widetilde{B},
\widetilde{N}_{\C}, \widetilde{S})$ is a topological Tits system.

\bigskip
\noindent
{\it The Affine Weyl Group}.  As a Coxeter system, the affine Weyl group
has generating set $\stil =S \cup \{s_0\}$, where $s_0$ is reflection
across the hyperplane $\alpha _0 =1$.  The cosets $\wtil /W$ have two
canonical sets of representatives: The coroot lattice $\cQ ^\vee$ and
the minimal length elements \wtils . An element $\lambda$ of the coroot
lattice is said to be {\it antidominant} if it lies in the closure of
$-\mathcal{C}$, the {\it antidominant} Weyl chamber.  We write
$\cQ^{\vee}_{-}$ for the set of antidominant elements in $\cQ^{\vee}$.  The
antidominant elements of the coroot lattice are characterized by the
following proposition.

\begin{proposition} \label{AntidominantCoroot}
For $\lambda \in \cQ^{\vee}$, the following conditions are equivalent:
\begin{enumerate}
\renewcommand{\labelenumi}{(\alph{enumi})}
\item $\lambda \in \wtils $.
\item $\lambda W$ is maximal in $W \cdot \lambda W/W$.
\item $\lambda$ is antidominant.
\end{enumerate}
\end{proposition}

\begin{proposition}
If $\sigma \in \widetilde{W}^S$ and $\lambda \in Q^{\vee}_{-}$,
then
$$
l(\sigma \lambda) = l(\sigma) + l(\lambda).
$$
\end{proposition}

%
%

\bigskip
\noindent
{\it The Affine Grassmannian.}  The affine Grassmannian associated to
$G$ is the homogeneous space $\affg = \gctil /P$. By the Iwasawa
decomposition, the natural map $\Omega _{alg} G \lra \affg$ is a \homeo
, and hence \affg\ is homotopy equivalent to $\Omega G$ by the
Quillen/Garland-Raghunathan theorem.  As the notation suggests, $\gctil
/P$ is closely analogous to a flag variety of a linear algebraic group.
For example, the $\widetilde{B}$-orbits in $\affg$ are
finite-dimensional affine cells called Schubert cells.  Each such cell
is of the form $e_{\sigma} = \widetilde{B} \cdot \sigma P/P$ where
$\sigma$ is an element of $\wtils$ or \coroot . More precisely, let
$P^- \subset \gctil$ denote the subgroup of maps regular at infinity,
let $P^{(-1)} =\{ f \in P^- : f(\infty) =I\}$, and let $\tilde{U}
_\sigma =\{u \in \util : \sigma ^{-1} u \sigma \in P^{(-1)}\}$. Then
$\tilde{U} _\sigma$ is a \fd\ unipotent group and hence is isomorphic to
an affine space; furthermore, the natural map $\tilde{U} _\sigma \sigma
\lra e_\sigma$ is an \iso\ of varieties. 

Given $\sigma \in \wtils$, if $n_{\sigma}$ is a lift of $\sigma$ to
$\widetilde{N}_{\C}$, we shall follow tradition and write $\sigma P$ to
denote the coset $n_{\sigma}P$.  Also, we will write $\sigma P/P$ for
the image of $\sigma P$ under the canonical projection $\gctil \lra
\affg$.

The closure $e_{\sigma}$, which is denoted $X_{\sigma}$, is called an
{\it affine Schubert variety} or simply, a Schubert variety.  These
varieties are (finite-dimensional) complex projective varieties.

\section{\affg\ as a Topological Group}

Since $\Omega _{alg} G$ is a topological group under pointwise
multiplication, \affg\ inherits a topological group structure from the
canonical \homeo\ $\Omega _{alg} G \stackrel{\cong}{\lra} \affg $. This
group structure is a special feature of our present, affine setting.
Indeed, a classical homogeneous space of the form $G_{\C}/P$ cannot even
be given the structure of an $H$-space.  A classical theorem of Hopf 
shows that the rational cohomology ring of a finite $H$-space is
isomorphic to an exterior algebra on odd-dimensional classes.  The
cohomology of $G_{\C}/P$, however, is concentrated in even dimensions.
                                                      
On the other hand, there is no algebraic ind-group structure on \affg
. Here {\it algebraic ind-group} means ``group object in the category of
ind-varieties''; it does not mean that the filtrations defining the
ind-variety structure are subgroups. We say that an ind-variety $X$ is
irreducible if it cannot be written as a union of two proper Zariski
closed subsets. It is easy to show that $X$ is irreducible if and only
if it admits a filtration by finite-dimensional irreducible
subvarieties.

\begin{proposition}
\label{AbelianInd-Variety}
Let $G$ be an irreducible projective algebraic ind-group.
Then $G$ is abelian.
\end{proposition}

\proof Let $G_n$ be a filtration of $G$ by irreducible projective
varieties (note that $G_n$ need not be a subgroup). Let $\gamma : G
\times G \lra G$ be the map that sends $(g,h)$ to its commutator
$ghg^{-1}h^{-1}$ and observe that $\gamma$ is a morphism of
ind-varieties.  For $h \in G$, we let $\gamma_h(g) = \gamma(g,h)$ and
note that $\gamma_e(g) = e$ for every $g \in G$.  Since each $G_n$ is
irreducible and projective, Lemma 1 of [Shafarevich, 4.3] applies to the
restriction of $\gamma$ to $(G \times G)_n = G_n \times G_n$.  Thus,
$\gamma_h$ is constant for each $h \in G_n$.  Since $\gamma_h(e) = e$,
it follows that each $\gamma_h \equiv e$, whence $\gamma(g,h) = e$ for
all $g, h$ in $G_n$. Since $n$ may be chosen arbitrarily, we conclude
that $G$ is abelian.

\begin{corollary}
The topological group structure on $\LG$ is not an algebraic ind-group
structure. 
\end{corollary}

We can say more, however.

\begin{proposition}
\label{RegularInversion}
Let $G$ be an ind-variety with filtration $(G_n)$.  Suppose that $G$
admits a topological group structure with a regular multiplication
map.  If each $G_n$ is a normal variety, then the inversion map is also
regular.
\end{proposition}

\proof Let $G$ satisfy the hypotheses stated in the proposition.  Let
$\mu : G \times G \lra G$ be the multiplication map and $\chi : G \lra
G$ be the inversion map.  For a fixed level $G_p$ of the filtration,
there exist integers $q$ and $n$ for which $\chi(G_p) \subset G_q$ and
$\mu:G_p \times G_q \lra G_n$ is a regular map of finite-dimensional
varieties.  Let $V \subset G_p \times G_q$ be the fiber of $\mu$ over
the identity element $e$.  Since $\mu$ is regular, $V$ is a subvariety
of $G_p \times G_q$.  The projection $\pi_1 : V \lra G_p$ is
bijective and, since $G_p$ is normal, an isomorphism (\cite{kumar},
p. 514).  We conclude that $\chi = \pi_2 \pi_1^{-1} : G_p \lra G_q$ is
algebraic.

\begin{corollary} \label{nonexist}
There does not exist a topological group structure on $\LG$ with a
regular multiplication map.
\end{corollary}

\proof Suppose that such a structure exists.  Since the levels of the
standard filtration on $\LG$ are Schubert varieties---and thus
normal (\cite{kumar}, p. 274)---Propositions \ref{RegularInversion}
and \ref{AbelianInd-Variety} imply that the group structure is
abelian.  It follows that $\LG$ is homotopy-equivalent to a product of
Eilenberg-MacLane spaces (\cite{hatcher}, 4K.7).  Since the cohomology of
$\LG$ is torsion-free and concentrated in even degrees, the
Eilenberg-MacLane spaces in question must be of type $K(\Z,2)$.  Thus
$\LG$ is a product of $\C P^{\infty}$'s.  Since $\pi_2(\LG) = \Z$,
$\LG \cong \C P^{\infty}$, a contradiction in any Lie type.

\bigskip

We next show that the inversion map on $\LG$ is not regular.  The
proof makes use of the following lemma.

\begin{lemma}
For $w \in W$, let $\varphi_w : G/T \lra G/T$ denote the map $gT
\mapsto gwT$.  If $w \neq 1$, then $\varphi_w$ is not a morphism of
varieties $\GC /B \lra \GC /B$.
\end{lemma}

\proof Suppose that $w \ne 1$ and that $\varphi_w : \GC /B \lra \GC
/B$ is regular.  Let $r_w : \GC \lra \GC$ denote right multiplication
by $w$ and consider the following diagram:
$$\begin{diagram}
\node{\GC} \arrow{e,t}{r_w} \arrow {s,l}{\pi}
\node{\GC} \arrow{s,r}{\pi} \\
\node{\GC /B} \arrow{e,b}{\varphi_w} \node{\GC /B}
\end{diagram}$$

\noindent
Observe that the diagram commutes when restricted to $G$ in the upper
left.  Since $G$ is Zariski dense in $G_{\C}$, the regularity of
$\varphi_w$ implies that the unrestricted diagram commutes.  This
yields a contradiction, however, since $\pi(U_w)$ is a point, whereas
$\pi \circ r_w$ is an embedding on $U_w$.

\begin{proposition}
The inversion map $\chi:\LG \lra \LG$ is not regular.
\end{proposition}

\proof Note that, for $\lambda \in Q^{\vee}$, $\chi$ takes the $P$-orbit
$\mathcal{O}_\lambda$ to the $P$-orbit $\mathcal{O}_{\lambda ^{-1}}$.
Fix $\lambda \in\cQ ^\vee _-$ with $\alpha _s \lam =0$ for all $s\in S$, and
let $Y_\lambda$ denote the $\GC$-orbit $G_{\C} \cdot \lambda$. The
corresponding subset of $\Omega _{alg}G$ is $Z_\lambda$, the image of
the map $\psi_\lambda :G/T \lra \Omega _{alg} G$ given by
$\psi_\lambda(gT) = g\lambda g^{-1}$.  Since $\lambda$ is dominant, the
dominant representative of the orbit $\mathcal{O}_{\lambda^{-1}}$ is
$w_0 \cdot \lambda ^{-1}$.  We thus have a commutative diagram of
homeomorphisms:

$$\begin{diagram}
\node{G/T} \arrow{e,t}{\psi_\lambda} \arrow {s,l}{\varphi_{w_0}}
\node{Z_\lambda} \arrow{s,r}{\chi} \\
\node{G/T} \arrow{e,b}{\psi_{w_0 \cdot \lambda^{-1}}} \node{Z_\lambda}
\end{diagram}$$

\bigskip
\noindent
Expressing this in terms of varieties yields 
$$\begin{diagram}
\node{\GC /B} \arrow{e,t}{} \arrow {s,l}{\varphi_{w_0}}
\node{Y_\lambda} \arrow{s,r}{\chi} \\
\node{\GC /B} \arrow{e,b}{} \node{Y_\lambda}
\end{diagram}$$

\bigskip
\noindent
where now the horizontal maps are just the usual orbit maps $gB
\mapsto g\lambda$, $gB \mapsto gw_0 \lambda ^{-1}$, respectively.  If
$\chi$ is algebraic in this diagram, then it is an isomorphism of
varieties.  This implies $\varphi_{w_0}$ is algebraic, contradicting
the lemma.

\section{The balanced product \balp }

In this section we study the balanced product $\cE _G =\balp $. We begin
with a discussion of the orbit map $\gctil \lra \affg $. Next we define
a natural ind-variety structure on \balp . Finally, we discuss the
Pontrjagin product and the balanced product in homology.

\subsection{The action map $\gctil \times \affg \lra \affg$}

The following proposition is a special case of \cite{kumar}, Lemma
7.4.10 and Proposition 7.4.11. 

\begin{proposition} \label{} \marginpar{}

The orbit map $\pi : \gctil \lra \affg$ is a morphism of
ind-varieties. Furthermore, $\pi$ is a Zariski locally trivial principal
$P$-bundle, as ind-varieties. 

\end{proposition}

Note that $\pi$ is not globally trivial as ind-varieties, since \gctil\
is affine and \affg\ is projective. On the other hand, it is a striking
fact that $\pi$ {\it is} globally trivial in the complex topology; this
follows immediately from the Iwasawa decomposition. In contrast, for a
finite-dimensional flag variety $\gc /Q_I$ the principal $Q_I$ bundle
$\gc \lra \gc /Q_I$ is never topologically trivial, since $\pi _2 \gc
\cong \pi _2 G=0$, whereas $\pi _2 (\gc /Q_I) \cong H_2 (\gc /Q_I) \cong
\bz ^{|S|-|I|}$. 

Call a morphism of ind-varieties $f :X \lra Y$ a {\it quotient morphism}
if whenever $Z$ is an ind-variety and $g: Y \lra Z$ is a set map such
that $gf$ is a morphism of ind-varieties, we have $g$ a morphism of
ind-varieties. Any regular map with local sections is a quotient
morphism, hence:

\begin{corollary} \label{} \marginpar{}

$\pi$ is a quotient morphism. 

\end{corollary}

This yields:

\begin{proposition} \label{} \marginpar{}

The action map $\theta: \gctil \times \affg \lra \affg$ is a morphism of
ind-varieties. 

\end{proposition}

The \gctil -action on \affg\ has two nice properties: First, if $Q
\subset \gctil$ is any proper, standard parabolic subgroup, then \affg\
is exhaustively filtered by $Q$-invariant \svars\ \xlam . Second, the
action of $Q$ on such an \xlam\ factors through a finite-dimensional
quotient algebraic group $Q/Q_\lambda$. In other words, the action is
``regular'' in the sense of \cite{kumar}, p. 213, except that our
parabolics are not pro-groups but rather subgroups of certain associated
pro-groups. 

In the case of interest to us, namely $Q=P$, there is a simple explicit
description of the finite-dimensional quotients and the factored
action. Thinking of $P$ as $G(\bc [z])$, there are natural surjective
morphisms of algebraic ind-groups $\epsilon _k : P \lra G[k]$, where
$G[k] =G(\bc [z]/z^k)$. Let $P^{(k)}$ denote the kernel of $\epsilon
_k$. Then if $\lambda \in \antidom$, it is not hard to show directly
that $P^{(k)}$ acts trivially on \xlam\ for all sufficiently large
$k$. Hence the finite-dimensional algebraic group $G[k]$ acts
algebraically on \xlam . 

We note also that $\gctil /P^{(k)} \lra \affg$ is a Zariski locally
trivial principal $G[k]$-bundle, as ind-varieties. Indeed for this it is
enough to show that $\epsilon _k$ is locally trivial; but in fact it is
globally trivial, as a principal $P^{(k)}$-bundle of ind-varieties
(\cite{gm}).

\subsection{Ind-variety structure on $\cE _G$}

As topological spaces, we may filter \balp\ by the subbundles $\gctil
\times _P \xlam \lra \affg$, where \lam\ ranges over \antidom
. Restricting these subbundles to $X_\sigma$, $\sigma \in \wtils$, we
obtain a local product $E_{\sigma , \lambda} \lra X_\sigma$ with fiber
\xlam . In fact \esiglam\ has a unique structure of projective variety,
such that $E_{\sigma , \lambda} \lra X_\sigma$ is locally trivial in the
Zariski topology as varieties. To prove this we need two lemmas. 

\begin{lemma} \label{} \marginpar{}

Let $\pi : E \lra X$ be a morphism of noetherian schemes over a noetherian
base scheme $A$, and suppose that $\pi $ is locally trivial with fiber
$F$. Let \cP\ be any class of morphisms over $A$ that is closed under
composition and base change, and local in the target. Then if $F \lra A$
and $X\lra A$ are in \cP , so is $E \lra A$.

In particular this holds when \cP\ is the class of separated, finite
type, or proper morphisms. 

\end{lemma}

\proof The first part is an easy exercise. The last assertion follows
from standard facts (\cite{hartshorne}, \S II.4). 

\begin{lemma} \label{} \marginpar{}

Let $\pi : E \lra X$ be a local product with fiber $F$ in the complex
topology, and suppose that (i) $X, F$ are varieties; (ii) there are
local trivializations on a Zariski open cover $U_\alpha$ of $X$; and
(iii) the transition functions $U_\alpha \cap U_\beta \times F \lra
U_\alpha \cap U_\beta \times F $ associated to the cover in (ii) are
regular. Then $E$ has a unique algebraic variety structure compatible
with the trivializations in (ii), so that $\pi$ is a local product as
varieties.  Moreover, if $X$ and $F$ are complete, then so is $E$.

\end{lemma} 

\proof Identify the category of varieties over \bc\ with the category of
reduced separated schemes of finite type over $Spec \, \bc$. By the
gluing lemma (\cite{hartshorne}, II.2.2) we obtain a reduced scheme
$E^\prime$ over $Spec \, \bc$, and a morphism of schemes $E^\prime \lra
X$ that is locally trivial with fiber $F$. By the previous lemma,
$E^\prime$ is separated of finite type, and is complete if $X$ and $F$
are complete. Furthermore, the set of complex points $E^\prime (\bc )$
is naturally identified with $E$. This proves the lemma.

\bigskip

This yields a complete ind-variety structure on \balp , indexed by the
poset $\wtils \times \cQ ^\vee _-$. Note that the maps $\pi , \phi: \balp \lra
\affg$ given by $\pi ([g,x]) =gP$, $\phi ([g,x])=gxP$ are morphisms
of ind-varieties. Then $\Psi =(\pi, \phi)$ defines an isomorphism of
ind-varieties 

$$\Psi : \balp \stackrel{\cong}{\lra} \affg \times \affg .$$

\noindent The inverse map is given by $(gP, xP) \mapsto [gP, g^{-1}
xP]$. Since \esiglam\ is complete, it maps isomorphically onto its image
in the projective ind-variety $\affg \times \affg$. Hence $\esiglam$ is
projective, and \balp\ is a projective ind-variety.

\bigskip

\noindent {\it Remarks:} 1. We could have used $\Psi$ to define the
ind-variety structure, but the direct construction given above is more
natural and has more general applicability. Note that $\Psi$ embeds
\esiglam\ as a subvariety of $X_\sigma \times X_{\sigma \lambda}$. 

2. A similar argument shows that the $G[k]$ bundle
$\gctil /P^{(k)} \lra \affg$ has a natural ind-variety structure, such
that if $P^{(k)}$ acts trivially on \xlam , then $\esiglam \lra \xlam$ is an
algebraic fiber bundle with structure group $G[k]$.

\subsection{Bi-Schubert cells in \esiglam }

Fix $\sigma \in \wtils$ and $\lambda \in \antidom$. Then \esiglam\ is
partitioned into {\it bi-Schubert cells}:

$$\esiglam =\coprod _{\tau \leq \sigma, \nu \leq \lambda} e_{\tau ,
  \nu},$$ 

\noindent where $e_{\tau ,\nu }=U_\tau \tau \times e _\nu$, identifying
the latter with a subspace of \esiglam\ in the evident way. 

\begin{lemma} \label{} \marginpar{}

a) Each bi-Schubert cell is a locally closed subvariety of \esiglam .

b) The closure relations on the bi-Schubert cells are given by the
product order on \linebreak $\wtils \times \wtils $. 

\end{lemma}

\proof The proof is similar to one of the standard proofs for the
Schubert cells, so we only sketch the method. First we recall the 
Bott-Samelson varieties (see \cite{kumar}, \S 7.1). Let
$\underline{s}=(s_1,...,s_p)$ be a sequence of elements of \stil , and
let $Z_{\underline{s}}$ denote the balanced product

$$Z_{\underline{s}} =P_{s_1} \times _{\tilde{B}} P_{s_2} \times
_{\tilde{B}} ...  P_{s_p} \times _{\tilde{B}} *,$$

\noindent where $P_{s_i}$ denotes the parabolic subgroup generated by
\btil\ and $s_i$, and $*$ is a point. The Bott-Samelson variety
$Z_{\underline{s}}$ is a smooth projective variety and an iterated $\bp
^1$-bundle. We write $[g_1,...,g_p]$ for an equivalence class defining a
point in $Z_{\underline{s}}$. Then there is a regular map 
$\rho : Z_{\underline{s}} \lra \affg $ given by $\rho [g_1,...,g_p] =g_1
g_2...g_p P$. If $w=s_1...s_p$ is an $S$-reduced expression, then $\rho$
is a smooth resolution of $X_w$; in fact it is an isomorphism over the
top cell $e_w$. 

Now suppose $\tau, \nu \in \wtils$ as above. If $w=\tau \nu$ is a
reduced product, choosing an $S$-reduced expression compatible with the
product yields a Bott-Samelson map $\tilde{\rho}: Z_{\underline{s}} \lra 
\overline{e}_{\tau, \nu}$ commuting in the diagram

\bigskip

$
\begin{diagram}
\node{Z_{\underline{s}}} \arrow{e,t}{\rho} \arrow{s,l}{\tilde{\rho}} 
\node{X_w}\\
\node{\overline{e_{\tau, \nu}}} \arrow{ne,b}{\phi}
\end{diagram}
$

\bigskip

The proof is completed by analyzing $\tilde{\rho}$; details are left to
the reader.

\subsection{The balanced product in homology}

We define the {\it balanced product} on $H_* \affg$ by 

$$[X_\tau] \star [X _\nu] =\phi _* [E_{\tau, \nu}].$$

\begin{proposition} \label{} \marginpar{}

The balanced product is associative but not commutative, and satisfies 

\[ [X_\tau] \star [X _\nu] =   \left\{ \begin{array}{ll}
[X_{\tau \nu}]    & \mbox{if $\tau \nu$ is $S$-reduced}\\
 0   & \mbox{otherwise}
\end{array}
\right. \]

\end{proposition}

\proof Suppose $\tau \nu$ is reduced. Then since $E_{\tau, \nu} \lra
X_{\tau \nu}$ is an \iso\ of varieties over the top cell, it has degree
one in homology. Hence $[X_{\tau}] \star [X_\nu] =[X_{\tau \nu}]$. If
$\tau \nu$ is not reduced, we still have a Bott-Samelson map
$\tilde{\rho} : Z \lra E_{\tau, \nu}$, so that $\rho _* [Z] =[E_{\tau
,\nu}]$. Then $\phi \circ \tilde{\rho} =\rho : Z \lra \affg $. Since the
product is not reduced, it follows easily that $dim \,\rho (Z) < dim \,
Z$ and hence $\rho _* [Z]=0$. This proves the displayed formula. Clearly
the product is associative, and is not commutative since $\tau \nu$
$S$-reduced need not imply $\nu \tau$ $S$-reduced.

\bigskip

\begin{proposition} \label{balancefacts} \marginpar{}                                       
Suppose $\sigma \in \wtils$, $\lambda \in \antidom$. Then                       
a) $\phi (\esiglam) =X_{\sigma \lambda}$                               

b) $\phi_* ([\esiglam ]) =[X_{\sigma \lambda}]$.

c) $\phi _* : H_* \esiglam \lra H_* X_{\sigma \lambda}$ is
onto. 

d) For every $\omega \leq \sigma \lam$, there exists $\tau \leq \sigma$
and $\nu \leq \lam$ in \wtils\ such that $[X_\tau] \star
[X_\nu]=[X_\omega]$. 
                    
\end{proposition}                                                             

\proof Parts a, b and c follow easily from the cell decomposition and
the Steinberg Lemma (\cite{steinberg}, Theorem 15, and \cite{kp},
Proposition 3.1). Part (d) is a refinement of (c), proved as follows:
Suppose $\omega \leq \sigma \lambda$, and choose an $S$-reduced product
decomposition $\omega =\tau \nu$ with $\tau \leq \sigma$ and $\nu \leq
\lambda$. Then $\nu \in \wtils$. We claim that if we choose such a
decomposition with $\nu$ of maximal length, then also $\tau \in \wtils$:
For if not, then $\tau \da \tau s$ for some $s \in S$, and hence there
is a reduced product $\tau =\tau ^\prime s$. But then $\omega =\tau
^\prime s \nu$ is reduced, and $s\nu \leq s \lam \leq \lam$, since $\lam
\in \antidom$. This contradicts the maximality of $\nu$, proving our
claim. Thus $[X_\tau] \star [X_\nu] =[X_\omega]$.

\subsection{The Pontrjagin product}

\begin{proposition} \label{toptriv} \marginpar{}

$\esiglam \lra X_\sigma$ is topologically a product bundle $X_\sigma
  \times \xlam$. In fact there is a unique trivialization $\eta$ such that the
  diagram 

\bigskip

$
\begin{diagram}
\node{X_\sigma \times X_ \lambda} \arrow{s,l}{m} \arrow{e,t}{\eta}
\node{\esiglam} \arrow{sw,r}{\phi}\\
\node{\affg}
\end{diagram}
$

\bigskip

\noindent commutes.

\end{proposition}

\proof This follows from the Iwasawa decomposition. 

\bigskip

\noindent {\it Remark:} These trivializations are not regular maps,
since if they were then the product map $m$ would be algebraic. We
conjecture that $\esiglam \lra X_\sigma$ is never a product bundle as
varieties, nor indeed a product at all.  We note also that $\eta$ does
not map $e_\sigma \times \elam$ into $e_{\sigma, \lambda}$; we only have
that $\eta (e_\sigma \times \elam) \subset U_\sigma \sigma \times \cO
_\lambda$.

\bigskip

\begin{proposition} \label{multkey} \marginpar{}

Suppose $\sigma \in \wtils$, $\lambda \in \antidom$. Then 

a) $m(X_\sigma \times \xlam) =X_{\sigma \lambda}$

b) $m_* ([X_\sigma ] \otimes [\xlam]) =[X_{\sigma \lambda}]$. 

\end{proposition}

\proof Part (a) follows from Proposition~\ref{balancefacts} and
Proposition~\ref{toptriv}. Since $\eta$ is a map of \xlam -bundles over
$X_\sigma$, having degree one on the base and the fiber, $\eta$ itself
has degree one. Part (b) then follows similarly. 

\section{The canonical generating variety}

A {\it generating complex} for an $H$-space $Y$ consists of a CW-complex
$K$ and a map $f: K \lra Y$ such that the image of $H_*f$ generates $H_*
Y$ as a ring. A {\it generating variety} for \affg\ is a
finite-dimensional projective subvariety $X$ such that $X \subset \affg$
is a generating complex. A \svar\ \xlam\ is a {\it geometric
generating variety} if (1) for all $n\geq 1$ the image of the $n$-fold
Pontrjagin multiplication $m^n :X_\lambda ^n \lra \affg$ is a \svar\
$X_{\lambda _n}$ with $dim \, X_{\lambda _n} =n \cdot dim \, \xlam$; and
(2) $m^n _* : H_* X_\lambda ^n \lra H_* X_{\lambda _n}$ is surjective.
Since the sequence $\lambda _n$ is cofinal in \wtils , a \ggenv\ is in
particular a Schubert \genv . Now let $\lambda _0 =-\alpha _0 ^\vee \in
\antidom$. The following theorem provides a canonical \ggenv\ in every
Lie type.

\begin{theorem} \label{} \marginpar{}

$X_{\lambda _0}$ is a geometric \genv\ for \affg . 

\end{theorem}

\proof Take $\lambda _n =-n \alpha _0 ^\vee$. The theorem then follows
from Proposition~\ref{multkey} and induction on $n$. 

\bigskip

As a corollary we obtain a classical theorem of Bott \cite{bott}. 

\begin{corollary} \label{} \marginpar{}

The homology ring $H_* \Omega G$ is \fg . 

\end{corollary}

Bott's proof uses \genc es; indeed he introduced the concept partly for
this purpose. His \genc es, however, are smooth manifolds, whereas ours
are always singular varieties as we will see shortly. We will also see
that smooth Schubert generating varieties exist if and only if $G$ is
not of type $E_8$, $F_4$, or $G_2$. Our canonical \genvs\ have a
number of interesting properties:

\bigskip

\noindent \xlamz\ {\it as Thom space}: Taking $P$-orbits instead of
\btil -orbits in \affg\ yields a stratification 

$$\affg =\coprod _{\lambda \in \cQ ^\vee _-} \cO _\lambda,$$

\noindent where $\cO _\lambda$ is isomorphic to a vector bundle $\xi
_\lambda$ over the {\it Levi orbit} $M_\lambda =\gc \cdot \lam P$. Here
$M_\lambda$ is a flag variety $\gc /Q_{I(\lambda)}$, where $I(\lambda )
=\{s \in S: \, s \lam W=\lam W\}$. For details of this construction, see
\cite{mpar}. 

Taking $\lam =\lambda _0$, one finds that $I(\lambda _0)$ corresponds to
the nodes on the ordinary Dynkin diagram that are not adjacent to $s_0$
in the affine diagram. The adjoint action of $Q=Q_{I(\lambda _0)}$ on
$\liegc $ stabilizes the one-dimensional root subalgebra $\frak{u}
_{\alpha _0}$, so that  

$$\xi _{\lambda _0} =\gc \times _Q \frak{u} _{\alpha _0}.$$

\noindent Furthermore, since $\lambda _0$ is the minimal nonzero
element of \antidom , we have $\xlamz =\overline{\cO} _{\lambda _0}=
\cO _{\lambda _0} \cup *$, where $*$ is the basepoint. Hence \xlamz\ is
the one-point compactification of $\cO _{\lambda _0}$ and so can be
identified with the Thom space $T(\xi _{\lambda _0})$. 

\bigskip

\noindent {\it Minimal nilpotent orbits.} Let $Y_0 \subset \liegc$ denote
the minimal non-trivial nilpotent orbit; that is, the orbit of a nonzero
element $x \in \frak{u} _{\alpha _0}$. Then 

$$Y_0 \cong \gc \times _Q (\frak{u} _{\alpha _0} -\{0\}).$$

\noindent In other words, $Y_0$ is the $\bc ^\times$ bundle associated
to $\xi _{\lambda _0}$. Hence $\xlamz -(M_{\lambda _0} \cup *) \cong
Y_0$.  Similarly, $\overline{Y}_0 \cong \xlamz -M_{\lambda _0} .$ In
particular, $\overline{Y}_0 $ and \xlamz\ have equivalent singular
points. Since the cone point in $\overline{Y}_0$ is known to be
singular, this gives one way to see that \xlamz\ is singular at the
basepoint. We will give a different proof below. For an elegant
computation of $H^* Y_0$, see \cite{juteau}.
  
\bigskip

\noindent {\it Birkhoff strata.} The Birkhoff stratification of \affg\
is ``dual'' to the Schubert cells, and defines a descending filtration
by ind-subvarieties $Z_\lambda$ ($\lambda \in \coroot \cong \wtils $)
\cite{kumar} \cite{ps} \cite{gm}. There is a unique codimension one
Birkhoff variety, namely $Z_{s_0} =Z_{\alpha _0 ^\vee}$. It is
$P$-invariant and hence $Z_{\alpha _0 ^\vee} \cap X_{\lambda _0} \subset
M_{\lambda _0}$; since both varieties are irreducible of codimension
one, it follows that the inclusion is an equality. Thus if $\cU _0$ is
the Zariski open complement $\affg -Z_{\alpha _0 ^\vee}$, we have 
$\cU _0 \cap \xlamz =\overline{Y}_0$. This observation enters into the
beautiful classification of minimal degeneration singularities in
\cite{mov}. 

\bigskip

\noindent {\it Segments.} Every $w \in \wtils$ has a canonical reduced
factorization $w=\sigma _1 \sigma _2 ...\sigma _n$ into {\it segments},
defined in \cite{bm2}. For present purposes we will use a slightly
different definition, ignoring the refinement for types $BD$ of
\cite{bm2}. We say that $\sigma \in \wtils$ is a {\it segment} if it is
in the left $W$-orbit of $s_0 $ in $\wtils \cong \wtil /W$. In other
words, the segments are precisely the elements that index the cells of
\xlamz , excluding the basepoint. Each such $\sigma$ has a unique
reduced factorization $\sigma =\nu s_0$ with $\nu \in W^J$, where $J$ is
the complement in $S$ of the set of nodes adjacent to $s_0$ in \stil
. It is then easy to show that each $w \in \wtils $ has a unique
factorization into segments $w=\sigma _1 \sigma _2 ...\sigma _n$ such
that each partial product $\sigma _1 ...\sigma _k$ is also in \wtils . 
The canonical generating variety can be viewed as a geometric expression
of this factorization; for example, it follows that there is a
corresponding unique $\star$-factorization of homology classes

$$[X_w] =[X_{\sigma _1}] \star [X_{\sigma _2}] \star ...\star [X_{\sigma
    _n}].$$ 

\noindent No such factorization holds for the Pontrjagin
product, however \cite{llms} \cite{lss}.  

\bigskip

\noindent {\it Examples:} (i) Type $A_n$, $n \geq 2$: This is the only
type in which $s_0$ is adjacent to more than one node in the affine
Dynkin diagram, and hence is the only type in which the flag variety
$M_{\lambda _0}$ is of non-maximal type. Explicitly, $M_{\lambda _0}$ is
the variety of flags of type $V^1 \subset V^n \subset \bc ^{n+1}$. Now
let $\xi _1, \xi _{n+1}$ denote the line bundles corresponding to $V^1$
and $\bc ^{n+1} /V^n$ respectively. Then $\xi _{\lambda _0} =Hom \, (\xi
_{n+1}, \xi _1 )$. In this case the canonical \genv\ is not very
efficient; the smallest \genvs\ are the minuscule \genvs\ $\bp ^n
\subset \affg$. 

(ii) Type $C_n$, $n \geq 1$ (type $A_1$ is best thought of as type $C_1$
in this context). Here $M_{\lambda _0}$ has type $C_n/C_{n-1}$; thus
$M_{\lambda _0} \cong \bp ^{2n-1}$. Then $\xi _{\lambda _0} \cong \eta
^{* 2}$ where $\eta$ is the canonical line bundle (in topologist's but
not geometer's terminology, i.e., $\eta ^*$ is the hyperplane section
bundle). Hence $c_1 (\xi _{\lambda _0}) \in H^2 \bp ^{2n-1}$ is twice a
generator. For later reference, we note that it follows that $X_{\lambda
_0}$ satisfies \pd\ rationally but not integrally.

Now let $\mu \in \coroot $ be the unique element immediately below
$\lambda _0=-\alpha _0 ^\vee$; that is, $\mu =s_1 \lambda _0$. Since $H_*
\Omega Sp(n) \cong \bz [a_1,...,a_n]$ with $|a_i|=4i-2$, it follows that
$X_\mu$ is a Schubert \genv\ for $\cL _{Sp(n)}$. In fact $X _\mu =T(\xi
_{\lambda _0} \downarrow {\Bbb C} P^{2n-2})$. This is the generating
complex of \cite{hopkins}. 

(iii) Type $G_2$. The flag variety $M_{\lambda _0}$ is $G_{2 {\Bbb
C}}/Q$, where $Q$ is the parabolic omitting the long node $s_2$. It is a
chain with cells indexed by $1, s_2, s_1 s_2, s_2 s_1 s_2, s_1 s_2 s_1
s_2, s_2 s_1 s_2 s_1 s_2$. For short, let $X_k$ denote the \svar\ of
dimension $k$ in $\gc /Q_{s_0}$, and define $y_k \in H^{2k} (\gc
/Q_{s_0})$ by $\langle y_k, [X_k]\rangle =1$. Then the Chevalley formula
shows that $y_1 =c_1 (\xi _{\lambda _0})$ and the cup products $y_1
y_{k-1} =a_k y_k$ with $a_k =1,3,2,3,1$ for $1 \leq k\leq
5$. Proposition~\ref{pd} below then shows that $X_{\lambda _0}$
satisfies \pd\ rationally but not integrally, and in particular is not
smooth.

\section{Smooth vs. singular generating varieties}

We first give an easy topological proof that \xlamz\ is singular. From
the Thom space description we see that \xlamz\ is smooth away from the
basepoint (i.e., the point at infinity). 

\begin{proposition} \label{} \marginpar{}

a) \xlamz\ never satisfies \pd\ over \bz . Hence the basepoint is a
singular point. 

b) \xlamz\ satisfies \pd\ over \bq\ if and only if $G$ has type $A_1$,
$C_n$, or $G_2$. 

\end{proposition}

This is a special case of \cite{bm1}, but it is worthwhile to sketch a
proof here. Let $X$ be a space satisfying \pd\ of dimension $2n$, with
$H^{odd} X=0$. We say that a space is {\it palindromic} if it satisfies
\pd\ additively, and is a {\it chain} if has \ppoly\ $1 +t ^2 +t^4
+...+t^{2n}$. The proof of the following lemma is elementary, and left
to the reader.

\begin{lemma} \label{pd} \marginpar{}

Let $\xi \da X$ be a complex line bundle over $X$. Then 

a) $T(\xi)$ is palindromic if and only $X$ is a chain if and only if
$T(\xi)$ is a chain.

b) $T(\xi)$ satisfies \pd\ if and only if $H^* X \cong \bz [c_1
  \xi]/(c_1 \xi) ^{n+1}.$

\end{lemma}

We also have:

\begin{lemma} \label{} \marginpar{}

$\gc/ Q_{I(\lambda _0)}$ is a chain if and only if $G$ has type $A_1$,
  $C_n$ or $G_2$. 

\end{lemma}

The proof is an easy and amusing exercise with Dynkin
diagrams. Consider, for example, type $E_8$: 

 \vspace{.3in}

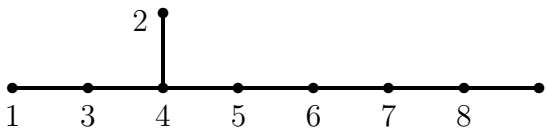
\raisebox{1ex}{\begin{picture}(5,1.2)
\mp(0,0)(1,0){8}{\ci}\put(2,1){\ci}
\put(0,0){\num{1}}\put(1,0){\num{3}}\put(2,0){\num{4}}\put(3,0){\num{5}}
\put(4,0){\num{6}}\put(5,0){\num{7}}\put(6,0){\num{8}} 
\put(2,.9){\makebox(0,0)[r]{$2$\hspace{.2\unitlength}}}
\put(2,0){\line(0,1){1}}
\put(0,0){\line(1,0){7}}
\end{picture}}
\\ \\ \vspace{.3in}

\noindent The Weyl group elements indexing the cells of $\gc/
Q_{I(\lambda _0)}$ begin with $s_8$, $s_7 s_8$, and so on until the fork
at $s_4$ is reached and we find two elements of length 6. Hence $\gc/
Q_{I(\lambda _0)}$ is not a chain in type $E_8$. Other types can be
ruled out in the same way; in the case of $F_4$ note that the node $s_3$
counts as a fork since one can reverse direction along the doubled edge.

The two lemmas imply that if \xlamz\ even rationally satisfies \pd ,
then $G$ must have type $A_1$, $C_n$, or $G_2$. As shown in the previous
section, these types do satisfy \pd\ rationally, but not integrally. 

\bigskip

We now recall the classical construction of Bott \cite{bott}, which
yields certain smooth \genvs\ for \affg . These \genvs\ are not always
\svars , but they are always either \svars\ or Levi orbits. Bott works
with the adjoint form $G^{ad}$ of $G$, so that $Hom \, (S^1, T^{ad})$
can be identified with the coweight lattice $\cP ^\vee$. If $\lam \in
\cP ^\vee$, let $C^{ad} _\lambda$ denote the centralizer of \lam\ in
$G$. The {\it Bott map} 

$$B_\lambda : G^{ad} /C^{ad} _\lambda \lra \Omega _0 G^{ad}$$

\noindent is defined by $B_\lambda (gC^{ad} _\lambda) =\lambda g 
\lambda ^{-1} g^{-1}$. Here $\Omega _0$ denotes the basepoint component
of the loop space. In fact Bott used the inverse of this map, but
$H_* B_\lambda$ and $H_ * B_\lambda ^{-1}$ generate the same subring of
$H_* \Omega _0 G^{ad}$, so the distinction will not be important. Note
that (i) the canonical map $\Omega G \lra \Omega _0 G^{ad}$ is a \homeo
; and (ii) $B_\lambda$ is in fact a map into the {\it algebraic} loops
$\Omega _{alg, 0} G^{ad}$. 

The main result of \cite{bott} is that whenever \lam\ is a ``short
minimal circle'', $G^{ad} /C^{ad} _\lambda \lra \Omega _0 G^{ad}$
is a \genc . In particular, the fundamental coweight $\omega _s ^\vee$
associated to a long simple root $\alpha _s$ is a short minimal circle;
we will confine our attention to this special case. If $\lam =\omega _s ^\vee$
happens to lie in the coroot lattice $\cQ ^\vee$ (examples occur in
every type except $A$, $C$), then $B_\lambda$ can be viewed in an
evident way as a commutator map $G/C_\lambda \lra \Omega _{alg}
G$. Identifying $\Omega _{alg} G =L_{alg} G/G =\affg$, we may also view
$B_\lambda$ as the map $gC_\lambda \mapsto \lambda g \lambda ^{-1}
P$. Since we are assuming $\lam \in \coroot$, multiplication by \lam\ is
homotopic to the identity and hence this last map is in turn homotopic
to the map $gC_\lambda \mapsto g \lambda ^{-1}P$. In other words, Bott's
\genc\ corresponds to the \gc -Levi orbit $M_{-\lambda }$. This yields
the following corollary of Bott's work: 

\begin{proposition} \label{} \marginpar{}

Suppose $\alpha _s$ is a long root and $\omega _s ^\vee \in
\coroot$. Then the \gc -Levi orbit $M_{-\alpha _s ^\vee }$ is a smooth \genv\
for \affg . 

In particular, if $G$ is not of type $A_n$ or $C_n$, then $M_{-\alpha _0
  ^\vee}$ is a \genv . 

\end{proposition}

The second assertion follows because except in types $AC$, there is a
unique node $t \in S$ adjacent to $s_0$, linked to $s_0$ by a single
edge. Hence $\alpha _t$ is long and $\omega _t ^\vee =\alpha _0 ^\vee$. 

\bigskip

At the opposite extreme, we can take $s$ to be a minuscule node; that
is, a node of $S$ that is in the orbit of $s_0$ under the action of the
automorphism group of the affine Dynkin diagram. (Warning: This is a
possibly non-standard use of the term {\it minuscule}.) The fundamental
coweights corresponding to the minuscule nodes are not in \coroot , and
indeed form a complete set of coset representatives for $\cP ^\vee
/\coroot $. Since these nodes always correspond to long roots, again
Bott's construction yields generating varieties for $\Omega _0
G^{ad}$. The second author showed in \cite{mbottfilt} that these
``minuscule \genc es'' correspond to smooth \svars\ that are in fact
geometric \genvs . In types $E_8$, $F_4$ and $G_2$, however, there are
no minuscule nodes, so this construction does not apply. In view of the
preceeding proposition, we nevertheless obtain:

\begin{proposition} \label{} \marginpar{}

In all Lie types, \affg\ admits a smooth \genv . 

\end{proposition}

Finally, the following proposition ties up a loose end from
\cite{mbottfilt}.

\begin{proposition} \label{notsmooth} \marginpar{}

\affg\ admits a smooth {\it Schubert} \genv\ if and only if $G$ is not
of type $E_8$, $F_4$ or $G_2$.

\end{proposition}

\proof What remains to be shown is that in the three exceptional types
no smooth \svar\ generates the homology. In fact we will show that no
smooth \svar\ can even generate the rational homology. Let $e_1 \leq
e_2 ...\leq e_r$ denote the exponents of $W$, where $r =|S|$, and recall
that $H_* (\Omega G ;\bq) \cong \bq [a_1,...,a_r]$, where
$|a_i|=2e_i$. Hence if \xlam\ is a rational \genc , we must have $dim \,
\xlam \geq e_r$. In the three exceptional types $e_r =29, 11, 5$
respectively \cite{bourbaki}. The smooth \svars\ were classified in
\cite{bm1}; in particular there are only finitely many in each Lie
type. The maximal dimension $d$ of a smooth \svar\ \xlam\ is given by:

$E_8$: $d=14$ (\xlam\ is a quadric of type $D_8/D_7$)

$F_4$: $d=7$ (\xlam\ is a quadric of type $B_4/B_3$)

$G_2$: $d=2$ ($\xlam \cong \bp ^2$).

Hence none of these three types has a smooth Schubert \genv .

\end{document}